\documentclass[preprint,12pt]{elsarticle}

\usepackage{graphicx,color}
\usepackage{amsfonts}

\usepackage{amssymb}
\usepackage{amsthm}

\usepackage{amsmath, mathrsfs, bm}

\newtheorem{thm}{Theorem}
\newtheorem{lem}[thm]{Lemma}

\newtheorem{cor}[thm]{Corollary}
\newdefinition{defn}{Definition}
\newdefinition{rmk}{Remark}
\newdefinition{alg}{Algorithm}
\newdefinition{exmp}{Example}
\newproof{pf}{Proof}

\begin{document}

\begin{frontmatter}




\title{Ideal Projectors of Type Partial Derivative and Their Perturbations}
\tnotetext[fund]{This work was supported by National Natural Science Foundation of China (No. 11171133 and 11101185).}


\author[Li]{Zhe Li}
\author[ZD]{Shugong Zhang}
\author[ZD]{Tian Dong\corref{cor1}}
\cortext[cor1]{Corresponding author} \ead{dongtian.jlu@gmail.com}

\address[Li]{Key Laboratory of Mathematics Mechanization, AMSS, Beijing 100190, China}
\address[ZD]{School of Mathematics,
Key Lab. of Symbolic Computation and Knowledge Engineering
\textup{(}Ministry of Education\textup{)}, Jilin University,
Changchun 130012, China}

\begin{abstract}

In this paper, we verify Carl de Boor's conjecture on ideal projectors for real ideal projectors of type partial derivative by proving that there exists a positive $\eta\in \mathbb{R}$ such that a real ideal projector of type partial derivative $P$ is the pointwise limit of a sequence of Lagrange projectors which are perturbed from $P$ up to $\eta$ in magnitude. Furthermore, we present an algorithm for computing the value of such $\eta$ when the range of the Lagrange projectors is spanned by the Gr\"{o}bner \'{e}scalier of their kernels w.r.t. lexicographic order.
\end{abstract}

\begin{keyword}
Ideal projector \sep Ideal projector of type partial derivative \sep Carl de Boor's conjecture

\MSC 41A10 \sep 41A35 \sep 13P10

\end{keyword}

\end{frontmatter}


\section{Introduction}
Polynomial interpolation is to construct a polynomial $p$ belonging
to a finite-dimensional polynomial subspace from a set
of data that agrees with a given function $f$ at the data set.
Univariate polynomial interpolation has a well developed theory,
while the multivariate one is very problematic since  a multivariate
interpolation polynomial is determined not only by the cardinal but
also by the geometry of the data set, cf. \cite{dBo94, GS00:2}.

As an elegant form of multivariate approximation, ideal interpolation
provides a natural link between multivariate polynomial interpolation and
algebraic geometry\cite{She2009}.   The study of
ideal interpolation was initiated by Birkhoff \cite{Bir1979}
and continued by several authors \cite{GS00:2,
dBo2005, She2009, LZD2011}.

Actually, ideal interpolation is an \emph{ideal projector} on polynomial ring whose
kernel is an ideal. When the kernel of an ideal projector $P$ is the vanishing ideal of
certain finite nonempty set $\Xi$ in $\mathbb{R}^d$, $P$ is a \emph{Lagrange
projector} on $\mathbb{R}[\bm{x}]:=\mathbb{R}[x_1,\ldots,x_d]$, the
polynomial ring in $d$ variables over $\mathbb{R}$, which
provides the Lagrange interpolation on $\Xi$. Obviously, $P$ is
finite-dimensional since its range is a $\#\Xi$-dimensional subspace
of $\mathbb{R}[\bm{x}]$. Lagrange projectors are standard
examples of ideal projectors.

It is well-known that every univariate ideal projector is an
\emph{Hermite projector}, namely it is the pointwise limit
of a sequence of Lagrange projectors. This inspired Carl de Boor\cite{dBo2005} to
conjecture that every finite-dimensional linear operator on
$\mathbb{C}[\bm{x}]$ is an ideal projector if and only if it is Hermite.

However, Boris Shekhtman\cite{BS2006} disproved this conjecture when the dimension $d\geq 3$.
In the same paper,  Shekhtman also showed that the conjecture is
true for bivariate complex projectors with the help of Fogarty
Theorem (see \cite{Foga1968}). Later, using linear algebra tools
only,  de Boor and Shekhtman\cite{deBoorShe2008} reproved the same result. Specifically,
Shekhtman\cite{BS2008} completely analyzed the bivariate ideal projectors
which are onto the space of polynomials of degree less than $n$ over
real or complex field, and verified the conjecture in this
particular case.


Let $P$ be an ideal projector that only interpolates a function and its partial derivatives. Obviously, many classical multivariate interpolation projectors are examples of $P$ which has applications in many fields of mathematics and science,  cf.\cite{Lor1992}. Naturally, we wonder whether de Boor' s conjecture is true for $P$ or not.


In this paper, a positive answer is offered to this question by Theorem \ref{mianthm} of Section \ref{mainresult} which states that there exists a positive $\eta\in \mathbb{R}$ such that $P$ is the pointwise limit of a sequence of Lagrange projectors which are perturbed from $P$ up to $\eta$ in magnitude, and the proof of the theorem is postponed to Section \ref{proof}, the last section of the paper. A further natural question is how to determine the value of $\eta$. We propose an algorithm in Section \ref{mainresult} for computing the value of such $\eta$ when the range of the Lagrange projectors is spanned by the Gr\"{o}bner \'{e}scalier of their kernels w.r.t. lexicographic order. And then, Section 4 is dedicated to some examples to illustrate the algorithm. The next section, Section 2, is devoted as a preparation for this paper.


\section{Preliminaries}\label{s:pre}

In this section, we will introduce some notation and review some
basic facts related to ideal projectors. For more details, we refer
the reader to \cite{dBo2005, She2009, deboor2006}.

Throughout the paper, we use $\mathbb{N}_0$ to stand for the monoid of
nonnegative integers and boldface type for tuples with their entries
denoted by the same letter with subscripts, for example,
$\bm{\alpha}=(\alpha_1,\ldots, \alpha_d)$.

Henceforward, we use $\leq$ to denote the usual product order on
$\mathbb{N}_0^d$, that is,
for arbitrary
$\bm{\alpha}$, $\bm{\beta}\in\mathbb{N}_0^d$, $\bm{\alpha}\leq \bm{\beta}$ if and only if $\alpha_i\leq \beta_i, i=1,\ldots, d$.
A finite nonempty set $\mathfrak{\Delta}\subset \mathbb{N}_0^d$ is called \emph{lower} if
for every $\bm{\alpha}\in \mathfrak{\Delta}$, $\bm{0}\leq\bm{\beta}\leq \bm{\alpha}$ implies  $\bm{\beta}\in \mathfrak{\Delta}$.

A \emph{monomial} ${\bm{x}}^{\bm{\alpha}}\in \mathbb{R}[\bm{x}]$ is
a power product of the form $x_1^{\alpha_1}\cdots x_d^{\alpha_d}$
with $\bm{\alpha}\in \mathbb{N}_0^d$. Thus, a \emph{polynomial} $p$ in
$\mathbb{R}[\bm{x}]$ can be expressed as a linear combination of monomials from $\mathrm{Supp}(p)$, the support of $p$, as follows,
\begin{equation}\label{pform}
p=\sum\limits_{\bm{\alpha}}\widehat{p}({\bm{\alpha}})
{\bm{x}}^{{\bm{\alpha}}}
\end{equation}
where $\widehat{p}({\bm{\alpha}})\in \mathbb{R}\backslash\{0\}$. For
$\bm{i}\in \mathbb{N}_0^d$ and $p \in \mathbb{R}[\bm{x}]$, if there
exists a monomial ${\bm{x}}^{\bm{\alpha'}}$ in $\mathrm{Supp}(p)$ such that $\bm{\alpha'} < {\bm{i}}$, then we
denote this fact as $p <_m \bm{i}$.

Let $P$ be a finite-dimensional ideal projector on $\mathbb{R}[\bm{x}]$.
The range and the kernel of $P$ are denoted by $\mathrm{ran}P$ and  $\mathrm{ker}P$ respectively.
Furthermore, $P$ has a dual
projector $P'$ on $\mathbb{R}'[\bm{x}]$, the algebraic dual of $\mathbb{R}[\bm{x}]$, whose
range can be described as
$$\mathrm{ran}P'=\{\lambda \in \mathbb{R}'[\bm{x}]: \mathrm{ker}P\subset \mathrm{ker}\lambda \},$$
which is the set of interpolation conditions matched by $P$.  Assume that $\Lambda\subset\mathbb{R}'[\bm{x}]$ is an
$\mathbb{R}$-basis for $\mathrm{ran}P'$,  then
$$\mathrm{ker}\Lambda:=\{f\in  \mathbb{R}[\bm{x}]: \lambda(f)=0, \forall\ \lambda\in \Lambda\}=\mathrm{ker}
P.$$

  We denote by
$\mathbb{T}^d$  the monoid of all monomials in $\mathbb{R}[\bm{x}]$.
For each fixed monomial order $\prec$ on $\mathbb{T}^d$,
a nonzero polynomial $f \in \mathbb{R}[{\bm{x}}]$ has a unique \emph{leading
monomial} $\mathrm{LM}_{\prec}(f)$, which is the $\prec$-greatest monomial appearing in $f$ with nonzero coefficient.
According to \cite{Mor2009}, the monomial set
$$
\mathcal{N}_\prec(\mathrm{ker}\Lambda):=\{{\bm{x}}^{\bm{\alpha}}\in \mathbb{T}^d:
\mathrm{LM}_{\prec}(f)\nmid {\bm{x}}^{\bm{\alpha}}, \forall f \in \mathrm{ker}\Lambda\}
$$
is
the \emph{Gr\"{o}bner \'{e}scalier} of $\mathrm{ker}\Lambda$ w.r.t. $\prec$.
We denote by
$\mathrm{ran}_{\prec}P$ the range of $P$ spanned by the Gr\"{o}bner \'{e}scalier of $\mathrm{ker}\Lambda$ w.r.t. $\prec$.

When $P$ is a Lagrange projector, we have $\mathrm{ker}\Lambda=\mathcal {I}(\Xi)$, the vanishing ideal of some finite nonempty set $\Xi\subset \mathbb{R}^d$.
In 1995, Cerlienco and Mureddu\cite{CM1995} proposed an purely combinatorial algorithm named MB for computing the Gr\"{o}bner \'{e}scalier of $\mathcal {I}(\Xi)$ w.r.t. some lexicographical order on $\mathbb{T}^d$ which is denoted by $\prec_{lex}$ here. Later, Felszeghy, R\'{a}th, and R\'{o}nyai\cite{FRR2006} provided a faster algorithm, lex game algorithm, by building a rooted tree $T(\Xi)$ of $d$ levels from $\Xi$ in the following way:
\begin{itemize}
  \item The nodes on each path from the root to a leaf
are labeled with the coordinates of a point.
  \item The root is regarded as the $0$-th level with no label,  its
children are labeled with the $d$-th coordinates of the points,
their children with the
$(d-1)$-coordinates, and so forth.
  \item If two points have same $k$
ending coordinates, then their corresponding paths coincide until
level $k$.
\end{itemize}
Given finite nonempty point sets $\Xi^{(1)}$, $\Xi^{(2)}\subset \mathbb{R}^d$ with
$\#\Xi^{(1)}=\#\Xi^{(2)}$. If $T(\Xi^{(1)})$ and $T(\Xi^{(2)})$ have same structure, \cite{FRR2006} showed that $\mathcal {N}_{\prec_{lex}}(\mathcal {I}(\Xi^{(1)}))=\mathcal {N}_{\prec_{lex}}(\mathcal {I}(\Xi^{(2)}))$.

\section{Main results}\label{mainresult}

Let
$$\delta_{\bm{\xi}}: \mathbb{R}[\bm{x}]\rightarrow \mathbb{R}: f\mapsto f(\bm{\xi})$$
denote the evaluation functional at the point
$\bm{\xi}=(\xi_1,\ldots,\xi_d)\in \mathbb{R}^d$, and let
$$\mathrm{D}^{\bm{\alpha}}: \mathbb{R}[\bm{x}]\rightarrow \mathbb{R}[\bm{x}]: f\mapsto \frac{\partial ^{\bm{\alpha}}}{\partial
{\bm{x}}^{\bm{\alpha}}}f:=\frac{\partial ^{\alpha_1+\cdots+\alpha_d}}{\partial
x_1^{\alpha_1}\cdots\partial x_d^{\alpha_d}}f$$ be the differential operator with
respect to $\bm{\alpha}=(\alpha_1,\ldots, \alpha_d)\in \mathbb{N}_0^d$ with $\mathrm{D}^{\bm{0}}=\mathrm{I}$, the identity operator on $\mathbb{R}[\bm{x}]$.


\begin{defn}\label{HermiteProjectorde}
Let $P$ be a finite-dimensional ideal projector on $\mathbb{R}[\bm{x}]$. If there exist distinct points
$\bm{\xi}^{(1)},\ldots,\bm{\xi}^{(\mu)}\in \mathbb{R}^d$ and their associated lower
sets
$\mathfrak{\Delta}^{(1)},\ldots,\mathfrak{\Delta}^{(\mu)}\subset
\mathbb{N}_0^d$ such that
\begin{equation}\label{HermiteProjector}
 \mathrm{ran}P'=\mathrm{Span}_{\mathbb{R}}\{\delta_{\bm{\xi}^{(k)}}\circ \mathrm{D}^{\bm{\alpha}}:  {\bm{\alpha}}\in \mathfrak{\Delta}^{(k)}, 1\leq
k\leq \mu\},
\end{equation}
namely $P$ only interpolates a function and its partial derivatives,
then we call $P$ an \emph{ideal projector of type partial
derivative}.
\end{defn}

As typical examples, Hermite projectors of type \emph{total degree} and of type \emph{coordinate degree} are both ideal projectors of type partial derivative, cf. \cite{Lor2000}.

\begin{lem}\label{delta0condition}
Let $\bm{\xi}^{(1)},\ldots,\bm{\xi}^{(\mu)}\in \mathbb{R}^d$ be
distinct points, and let
$\mathfrak{\Delta}^{(1)},\ldots,\mathfrak{\Delta}^{(\mu)}\subset \mathbb{N}_0^d$ be their associated lower
sets. Set
\begin{equation}\label{eta}
\eta_0:=\min\left\{\frac{\|\bm{\xi}^{(k)}-\bm{\xi}^{(l)}\|_2}{\|{\bm{\alpha}}-{\bm{\alpha}}'\|_2}:
\bm{\alpha} \in \mathfrak{\Delta}^{(k)}, \bm{\alpha}' \in
\mathfrak{\Delta}^{(l)},\bm{\alpha}\neq \bm{\alpha}', 1\leq k<l\leq \mu
\right\}.
\end{equation}
Then for arbitrary nonzero $h\in (-\eta_0, \eta_0)\subset \mathbb{R}$, the
point set
\begin{equation}\label{pointset}
\Xi_h:=\left\{\bm{\xi}^{(k)}+ h \bm{\alpha}: \bm{\alpha}\in
\mathfrak{\Delta}^{(k)}, 1\leq k\leq \mu\right\}
\end{equation}
exactly consists of
$\#\sum\limits_{i=1}^{\mu}\mathfrak{\Delta}^{(i)}$ distinct points.
\end{lem}

\begin{pf}
Suppose that there exist $\bm{\alpha}\in\mathfrak{\Delta}^{(k)}$ and $\bm{\alpha}'\in\mathfrak{\Delta}^{(l)}$
with $1\leq k<l\leq\mu$ such that $\bm{\xi}^{(k)}+ h\bm{\alpha}=\bm{\xi}^{(l)}+ h \bm{\alpha}'$ which implies that
$\bm{\alpha}\neq \bm{\alpha}'$ by $\bm{\xi}^{(k)}\neq \bm{\xi}^{(l)}$. Consequently, we have
$$h=\frac{\|\bm{\xi^{(k)}}-\bm{\xi^{(l)}}\|_2}{\|\bm{\alpha}-\bm{\alpha}'\|_2},
$$ which is in direct contradiction to the hypothesis that $0<|h|<\eta_0$. \qed
\end{pf}

Lemma \ref{delta0condition} holds out the possibility of
intuitively perturbing
an ideal projector of type
partial derivative to a sequence of Lagrange projectors.

\begin{defn}\label{LagrangeProjectorde}
Let $P$ be an ideal projector of type partial derivative on $\mathbb{R}[\bm{x}]$ with
$\mathrm{ran}P'$ described by (\ref{HermiteProjector}).
For an
arbitrary fixed $h\in \mathbb{R}$ with $0<|h|<\eta_0$ where $\eta_0$
is as in (\ref{eta}),  define $P_h$ to be the Lagrange projector on
$\mathbb{R}[\bm{x}]$ with
\begin{equation}\label{LagrangeProjector}
\mathrm{ran}
P'_h=\mathrm{Span}_{\mathbb{R}}\{\delta_{\bm{\xi}^{(k)}+ h \bm{\alpha}}:
\bm{\alpha}\in \mathfrak{\Delta}^{(k)},  1\leq k\leq \mu \}.
\end{equation}
Then $P_h$ is called an \emph{$h$-perturbed Lagrange projector of
$P$}.
\end{defn}

\begin{rmk}\label{dualbasis}
It is easy to see from (\ref{HermiteProjector}) and (\ref{LagrangeProjector}) that
$$\bm{\lambda}:=(\delta_{\bm{\xi}^{(k)}}\circ \mathrm{D}^{\bm{\alpha}}: \bm{\alpha}\in
\mathfrak{\Delta}^{(k)},~k=1,\ldots,\mu) \in
{(\mathbb{R}'[\bm{x}])}^n$$
and
$$\bm{\lambda}_h:=(\delta_{\bm{\xi}^{(k)}+ h
\bm{\alpha}}: \bm{\alpha}\in \mathfrak{\Delta}^{(k)}, k=1,\ldots,\mu) \in
{(\mathbb{R}'[\bm{x}])}^n$$
form $\mathbb{R}$-bases for $\mathrm{ran} P'$ and $\mathrm{ran} P'_h$ respectively, where $n=\sum_{k=1}^\mu \#\mathfrak{\Delta}^{(k)}$.

Moreover, an ordering $\prec_\lambda$ for the entries of $\bm{\lambda}$ and $\bm{\lambda}_h$ will be defined as follows: We say $\delta_{\bm{\xi}^{(k)}}\circ \mathrm{D}^{\bm{\alpha}}\prec_\lambda \delta_{\bm{\xi}^{(k')}}\circ \mathrm{D}^{\bm{\alpha'}}$ or $\delta_{\bm{\xi}^{(k)}+ h
\bm{\alpha}}\prec_\lambda \delta_{\bm{\xi}^{(k')}+ h\bm{\alpha'}}$ if
$$
k<k', \quad \mbox{or}\quad k=k' \mbox{ and } \bm{\alpha} \prec \bm{\alpha'},
$$
where $\prec$ is an arbitrary monomial order on $\mathbb{N}_0^d$.

\end{rmk}

We are now ready to give one of our main theorem, Theorem \ref{mianthm}, which states that every
ideal projector of type partial derivative on $\mathbb{R}[\bm{x}]$ is the pointwise limit of Lagrange projectors, namely Carl de Boor's conjecture is true for this type of ideal projectors.

\begin{thm}\label{mianthm}

Let $P$ be an ideal projector of type partial derivative on $\mathbb{R}[\bm{x}]$ with
$\mathrm{ran}P'$ described by \eqref{HermiteProjector}, and let \textup{(}$P_h$, $0<|h|<\eta_0$\textup{)} be a sequence of perturbed Lagrange projector of $P$ where $\eta_0$ is as in \textup{(\ref{eta})}. Then the following statements hold:
\begin{enumerate}
  \item[\textup{(i)}] There exists a positive $\eta\in \mathbb{R}$ such that
   $$
   \mathrm{ran}P_h=\mathrm{ran}P, \quad \forall 0<|h|<\eta\leq\eta_0.
   $$
  \item[\textup{(ii)}] $P$ is the pointwise limit of $P_h, 0<|h|<\eta,$ as $h$ tends to zero.
\end{enumerate}
\end{thm}

The proof of Theorem \ref{mianthm} will be provided in Section \ref{proof}. Actually, with similar methodology there, we can easily prove the following theorem, which is a more general version of Theorem
\ref{mianthm}.

\begin{thm}\label{corthm}
Let $P$ be an ideal projector of type partial derivative from $C^{\infty}(\mathbb{R}^d)$ onto $\mathrm{ran}P$, then
there exists Lagrange projector $P_h$ onto $\mathrm{ran}P$ such that
for all $f \in C^{\infty}(\mathbb{R}^d)$, $P f$ is the  limit of
$P_h f$  as $h$ tends to zero.
\end{thm}

Now, after introducing Definition \ref{borderbasisdef}, we have an immediate corollary of Theorem \ref{mianthm}.

\begin{defn}\label{borderbasisdef}\textup{\cite{She2009}}
Let $P$ be an ideal projector from $\mathbb{R}[\bm{x}]$
onto $\mathrm{ran}P$ with $\dim\mathrm{ran}P=n$. Assume that
$\bm{q}=(q_1,\ldots ,q_{n})\in {\mathbb{R}[\bm{x}]}^{n}$ is an
$\mathbb{R}$-basis for $\mathrm{ran} P$, and the border set
$\partial \bm{q}$ of $\bm{q}$ is defined by
$$\partial \bm{q}:=\{1,x_k q_l,k=1,\ldots,d,l=1,\ldots,n\}\setminus \{q_1,\ldots,q_{n}\}.$$
Then the set of polynomials
$$\{f-P f: f\in \partial \bm{q} \}$$
forms a \emph{border basis} for $\mathrm{ker}P$, which is
called a $\bm{q}$-\emph{border basis} for $\mathrm{ker}P.$
\end{defn}

\begin{cor}\label{miancor}
Let $P$ be an ideal projector of type partial derivative on $\mathbb{R}[\bm{x}]$, and let $\bm{q}$ be an $\mathbb{R}$-basis for $\mathrm{ran}P$. Then there exists a Lagrange projector $P_h$ onto $\mathrm{ran}P$  such that the $\bm{q}$-border basis for
$\mathrm{ker}P$ is the limit of $\bm{q}$-border basis for
$\mathrm{ker}P_h$ as $h$ tends to zero.
\end{cor}

%

Theorem \ref{mianthm} tells us that every ideal projector of
type partial derivative is the pointwise limit of Lagrange
projectors. Unfortunately, the converse statement is not true in general as the following example illustrates.

\begin{exmp}
Let $(P_h, 0<|h|<1)$ be a sequence of Lagrange projectors
with
\begin{align*}
\mathrm{ran} P_h&=\mathrm{Span}_{\mathbb{R}}\{1,x_1,x_2,x_1^2,x_1 x_2,x_2^2\},\\
\mathrm{ran} P'_h&=\mathrm{Span}_{\mathbb{R}}
\{\delta_{(0,0)}, \delta_{(0,h)}, \delta_{(h,0)}, \delta_{(1,1)}, \delta_{(1,1+h)}, \delta_{(1+h,1)}\},
\end{align*}
and let $P$ be an ideal projector with
\begin{align*}
\mathrm{ran} P'=\bigg\{&\delta_{(0,0)}\circ \mathrm{D}^{(0,0)},
\delta_{(0,0)}\circ \mathrm{D}^{(1,0)}, \delta_{(0,0)}\circ \mathrm{D}^{(0,1)},\\
&\delta_{(1,1)}\circ \mathrm{D}^{(0,0)}, \delta_{(1,1)}\circ \mathrm{D}^{(1,0)},
\delta_{(1,1)}\circ \mathrm{D}^{(0,1)}\bigg\}.
\end{align*}
However,  $\{1,x_1,x_2,x_1^2,x_1
x_2,x_2^2\}$ can not form an $\mathbb{R}$-basis for $\mathrm{ran} P$.
Hence, $(P_h, 0<|h|<1)$ can not converge pointwise to $P$, as $h$ tends to zero.
\end{exmp}

Consider the bijection
\begin{align*}
u:\mathbb{R}^d\times \mathbb{N}_0^d&\longrightarrow
{(\mathbb{R}\times\mathbb{N}_0)}^d\\
(\bm{\xi},\bm{\alpha})&\longrightarrow((\xi_1, \alpha_1),\ldots,(\xi_d, \alpha_d)).
\end{align*}
Let $\bm{\xi}^{(1)},\ldots,\bm{\xi}^{(\mu)}\in\mathbb{R}^d$ be distinct points and $\mathfrak{\Delta}^{(1)},\ldots,\mathfrak{\Delta}^{(\mu)}\subset
\mathbb{N}_0^d$ be lower sets. Then
\begin{equation}\label{Hermitetree}
\Omega:=\{u(\bm{\xi}^{(k)},\bm{\alpha}): \bm{\alpha}\in \mathfrak{\Delta}^{(k)},
k=1, \ldots, \mu\}\subset{(\mathbb{R}\times\mathbb{N}_0)}^d
\end{equation}
is called  an \emph{algebraic multiset}.  As mentioned by \cite{CM1995}, MB algorithm can be applied for
the algebraic multiset $\Omega$ to obtain
the Gr\"{o}bner \'{e}scalier of
the ideal
$$\{p\in \mathbb{R}[\bm{x}]: \delta_{\bm{\xi}^{(k)}} \circ \mathrm{D}^{\bm{\alpha}}(p)=0, \bm{\alpha}\in \mathfrak{\Delta}^{(k)}, \ 1\leq k\leq \mu\}$$
w.r.t. lexicographic order.

Recall Section \ref{s:pre}. We have known how to build a $d$-level tree $T(\Xi)$ from a finite nonempty set $\Xi\in \mathbb{R}^d$. If the space $\mathbb{R}^d$ is changed to ${(\mathbb{R}\times\mathbb{N}_0)}^d$, it is easy to see that we can also build a $d$-level tree $T(\Omega)$ from algebraic multiset $\Omega$ following the same rules, which makes lex game algorithm involved and leads to the following useful lemma.



\begin{lem}\label{mainlem}
Let $P$ be an ideal projector of type partial derivative with
$\mathrm{ran}P'$ as in \textup{(\ref{HermiteProjector})},  and let
$P_h$ be a perturbed
Lagrange projector of $P$. Let algebraic multiset
$\Omega\subset{(\mathbb{R}\times\mathbb{N}_0)}^d$ be as in \textup{(\ref{Hermitetree})} and  $\Xi_h\subset \mathbb{R}^d$ be as in \textup{(\ref{pointset})}.
If the rooted trees $T(\Omega)$ and $T(\Xi_h)$ have
the same structure, then
$$\mathrm{ran}_{\prec_{lex}}P=\mathrm{ran}_{\prec_{lex}} P_h.$$
\end{lem}

Next, we can proceed with another main theorem of this paper.

\begin{thm}\label{lexthm}
Let $P$ be an ideal projector of type partial derivative with
$\mathrm{ran}P'$ as in \textup{(\ref{HermiteProjector})},
and let $\textup{(}P_h, 0<|h|<\eta\textup{)}$ be a sequence
of $h$-perturbed Lagrange projectors of $P$, where $\eta$ is
obtained through Algorithm \textup{\ref{cond}} in the following.
If the range of  $P_h$ is $\mathrm{ran}_{\prec_{lex}}P_h$, then
the sequence $(P_h,
0<|h|<\eta )$
converges pointwise to the ideal projector $P$, as $h$ tends to zero.
\end{thm}

\begin{alg}\label{cond}(The range for $|h|$)
\vskip 3mm
\textbf{Input}: Distinct points
$\bm{\xi}^{(1)},\ldots,\bm{\xi}^{(\mu)}\in \mathbb{R}^d$  and lower
sets $\mathfrak{\Delta}^{(1)},\ldots,\mathfrak{\Delta}^{(\mu)}\subset \mathbb{N}_0^d$.\\
\indent\textbf{Output}:  A nonnegative number $\eta\in \mathbb{R}$ or $\infty$.\\
\indent\textbf{Step 1} Construct algebraic multiset $\Omega$ from $\bm{\xi}^{(1)},\ldots,\bm{\xi}^{(\mu)}$ and $\mathfrak{\Delta}^{(1)},\ldots,\mathfrak{\Delta}^{(\mu)}$ following (\ref{Hermitetree}), and then build rooted tree $T(\Omega)$ from $\Omega$ in the way introduced in Section \ref{s:pre}.\\
\indent\textbf{Step 2} Suppose that the first level nodes of $T(\Omega)$ are labeled with
the points of set $\mathcal {L}_{1}\subset \mathbb{R}\times\mathbb{N}_0$. \\
\indent\indent\textbf{Step 2.1} If $\#\mathcal {L}_{1}=1$, then $\eta\leftarrow \infty$.\\
\indent\indent\textbf{Step 2.2} If every point in $\mathcal {L}_1$ has the same first coordinate or the same second coordinate, then $\eta\leftarrow \infty$. \\
\indent\indent{\textbf{Step 2.3}}\textup{:} Otherwise, set
\begin{align*}
\eta\leftarrow\min\Bigg\{\frac{
|\xi_d^{(i)}-\xi_d^{(j)}|}{|\alpha_d^{(i)}-\alpha_d^{(j)}|}:&
\xi_d^{(i)}\neq\xi_d^{(j)}, \alpha_d^{(i)}\neq\alpha_d^{(j)}, \\ &(\xi_d^{(i)}, \alpha_d^{(i)})\mbox{ and }(\xi_d^{(j)}, \alpha_d^{(j)})\in \mathcal{L}_1\Bigg\}.
\end{align*}
\indent\textbf{\textup{Step 3}} Set $k\rightarrow 2$.\\
\indent\textbf{\textup{Step 4}} Suppose that the $k$-th level nodes are labeled respectively with the points of sets
$\mathcal {L}_{k}^{(1)}, \ldots, \mathcal {L}_{k}^{(\nu)}\subset\mathbb{R}\times\mathbb{N}_0$, where for each
$1\leq l\leq \nu$, the nodes labeled with the points in $\mathcal {L}_{k}^{(l)}$ share the same parent.
For $l=1, \ldots, \nu$ and $\#\mathcal {L}_{k}^{(l)}\geq 2$, do the following steps. \\
\indent\indent\textbf{\textup{Step 4.1}} Set
\begin{align*}
\eta'\leftarrow\min\Bigg\{\frac{
|\xi_{d-k+1}^{(i)}-\xi_{d-k+1}^{(j)}|}{|\alpha_{d-k+1}^{(i)}-\alpha_{d-k+1}^{(j)}|}:&
\xi_{d-k+1}^{(i)}\neq\xi_{d-k+1}^{(j)}, \alpha_{d-k+1}^{(i)}\neq\alpha_{d-k+1}^{(j)}\\
&(\xi_{d-k+1}^{(i)}, \alpha_{d-k+1}^{(i)})\mbox{ and }(\xi_{d-k+1}^{(j)}, \alpha_{d-k+1}^{(j)})\in \mathcal {L}_{k}^{(l)}\Bigg\}.
\end{align*}
\indent\indent\textbf{\textup{Step 4.2}} If $\eta'<\eta$, then $\eta\leftarrow \eta'$.\\
\indent\textbf{\textup{Step 5}} If $k=d$, then return $\eta$ and stop.
Otherwise set $k\leftarrow k+1$, continue with {\textup{Step 4}}.

\end{alg}



\begin{pf}

To prove this theorem,  by Lemma \ref{mainlem} and Theorem \ref{mianthm},
it suffices to
show that the rooted trees  $T(\Xi_h),  0<|h|<\eta$, and $T(\Omega)$ have
the same structure, where $\Xi_h$ is as in
(\ref{pointset}) and  $\Omega$ is as in (\ref{Hermitetree}).
Now, with the notation in Algorithm \textup{\ref{cond}}, we will use induction on the number
of levels $k$ of the rooted tree to prove this.

When $k=1$,  assume that
there exist some $(\xi^{(i)}_d, \alpha^{(i)}_d)$ and $(\xi_d^{(j)}, \alpha_d^{(j)})\in \mathcal {L}_1$ such that
$\xi_d^{(i)}+h \alpha_d^{(i)}= \xi_d^{(j)}+h
\alpha_d^{(j)}$. The same argument in Lemma \ref{delta0condition} shows that $h=
|\xi_d^{(i)}-\xi_d^{(j)}|/|\alpha_d^{(i)}-\alpha_d^{(j)}|$ where $\alpha_d^{(i)}\neq\alpha_d^{(j)}$ and $\xi_d^{(i)}\neq
\xi_d^{(j)}$, which contradicts
\begin{align*}
|h|<\min\Bigg\{\frac{
|\xi_d^{(i)}-\xi_d^{(j)}|}{|\alpha_d^{(i)}-\alpha_d^{(j)}|}:&
\xi_d^{(i)}\neq\xi_d^{(j)}, \alpha_d^{(i)}\neq\alpha_d^{(j)}, \\ &(\xi_d^{(i)}, \alpha_d^{(i)})\mbox{ and }(\xi_d^{(j)}, \alpha_d^{(j)})\in \mathcal{L}_1\Bigg\}.
\end{align*}
Hence, the first levels of
$T(\Xi_h),  0<|h|<\eta$, and $T(\Omega)$
have the same structure.

Suppose that the first $k-1$ levels of $T(\Xi_h),  0<|h|<\eta$, and $T(\Omega)$ have the same structure. Assume that there exists some
$1\leq l\leq \nu$ and $(\xi_{d-k+1}^{(i)}, \alpha_{d-k+1}^{(i)}), (\xi_{d-k+1}^{(j)}, \alpha_{d-k+1}^{(j)})\in\mathcal {L}_{k}^{(l)}$
such that
$\xi_{d-k+1}^{(i)}+h \alpha_{d-k+1}^{(i)}= \xi_{d-k+1}^{(j)}+h
\alpha_{d-k+1}^{(j)}$. Since $(\xi_{d-k+1}^{(i)}, \alpha_{d-k+1}^{(i)})$, $(\xi_{d-k+1}^{(j)},
\alpha_{d-k+1}^{(j)})$ have common parent, it is easy to see that $h=|\xi_{d-k+1}^{(i)}-\xi_{d-k+1}^{(j)}|/|\alpha_{d-k+1}^{(i)}-\alpha_{d-k+1}^{(j)}|$ where $\alpha_{d-k+1}^{(i)}\neq \alpha_{d-k+1}^{(j)}$ and $\xi_{d-k+1}^{(i)}\neq
\xi_{d-k+1}^{(j)}$, which contradicts the fact
\begin{align*}
|h|<\min\Bigg\{\frac{
|\xi_{d-k+1}^{(i)}-\xi_{d-k+1}^{(j)}|}{|\alpha_{d-k+1}^{(i)}-\alpha_{d-k+1}^{(j)}|}:&
\xi_{d-k+1}^{(i)}\neq\xi_{d-k+1}^{(j)}, \alpha_{d-k+1}^{(i)}\neq\alpha_{d-k+1}^{(j)}\\
&(\xi_{d-k+1}^{(i)}, \alpha_{d-k+1}^{(i)})\mbox{ and }(\xi_{d-k+1}^{(j)}, \alpha_{d-k+1}^{(j)})\in \mathcal {L}_{k}^{(l)}\Bigg\}.
\end{align*}
Therefore, the first $k$ levels of $T(\Xi_h),  0<|h|<\eta$, and $T(\Omega)$ have the same structure.
\qed
\end{pf}


\section{Example}\label{examples}

In this section, we will present several examples to illustrate Theorem \ref{lexthm}.

\begin{exmp}\label{ex1}
Assume that  $P_h$ is a Lagrange projector with
$$\mathrm{ran} P'_h=\mathrm{Span}_{\mathbb{R}}\{\delta_{(0,0)}, \delta_{(h,0)}, \delta_{(0,h)}, \delta_{(1,1)}, \delta_{(1+h,1)}, \delta_{(1,1+h)}\}.$$
Construct the rooted tree of the algebraic multiset
\begin{align*}
\Omega=\{&((0, 0), (0, 0)), ((0, 1), (0, 0)), ((0, 0), (0, 1)), \\
 &((1, 0), (1, 0)), ((1, 1), (1, 0)), ((1, 0), (1, 1))\}.
\end{align*}
$T(\Omega)$ is illustrated in Figure \textup{\ref{eg1}}.
\begin{figure}[!htbp]
\begin{center}
\includegraphics[width=8cm, height=5cm]{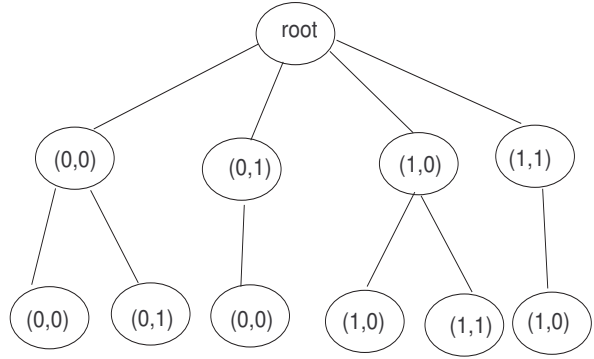}
\caption{$T(\Omega)$ of Example \ref{ex1}}\label{eg1}
\end{center}
\end{figure}
By Algorithm 1, we obtain $\eta=1$.
From Theorem \ref{lexthm}, we can conclude that
$(P_h, 0<|h|<1)$ onto $\mathrm{Span}_{\mathbb{R}}\{1, x_2, x_1, x_2^2,$\\$
x_1 x_2, x_2^3\}$ pointwise converges to an Hermite projector $P$ with
\begin{align*}
\mathrm{ran} P'=\{&\delta_{(0,0)}\circ \mathrm{D}^{(0,0)}, \delta_{(0,0)}\circ \mathrm{D}^{(1,0)}, \delta_{(0,0)}\circ \mathrm{D}^{(0,1)},
\delta_{(1,1)}\circ \mathrm{D}^{(0,0)}, \delta_{(1,1)}\circ \mathrm{D}^{(1,0)},\\
&\delta_{(1,1)}\circ \mathrm{D}^{(0,1)}\},
\end{align*}
as $h$ tends to zero.
\end{exmp}

\begin{exmp}\label{ex2}
Assume that  $P_h$ is a Lagrange projector with
\begin{align*}
\mathrm{ran} P'_h=\mathrm{Span}_{\mathbb{R}}\{&\delta_{(0,0,0)}, \delta_{(h,0,0)}, \delta_{(0,h,0)},\delta_{(0,0,h)},
\delta_{(1,1,1)}, \\
&\delta_{(1+h,1,1)}, \delta_{(1,1+h,1)},
\delta_{(1,1,1+h)}\}.
\end{align*}
Construct the rooted tree of the algebraic multiset
\begin{align*}
\Omega=\{&((0, 0), (0, 0), (0, 0)), ((0, 1), (0, 0), (0, 0)), ((0, 0), (0, 1), (0, 0)),\\
&((0, 0), (0, 0), (0, 1)), ((1, 0), (1, 0), (1, 0)), ((1, 1), (1, 0), (1, 0)), \\
 &((1, 0), (1, 1), (1, 0)), ((1, 0), (1, 0), (1, 1))\}.
\end{align*}
$T(\Omega)$ is illustrated in Figure \textup{\ref{eg2}}.
\begin{figure}[!htbp]
\begin{center}
\includegraphics[width=8cm, height=5cm]{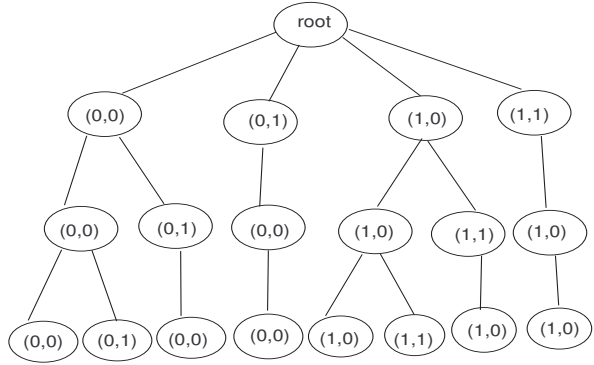}
\caption{$T(\Omega)$ of Example \ref{ex2}}\label{eg2}
\end{center}
\end{figure}
By Algorithm 1, we compute $\eta=1$.
From Theorem \ref{lexthm}, we can conclude that
$(P_h, 0<|h|<1)$ onto
$\{1, x_3, x_2, x_1, x_3^2, x_2
x_3,$\\$x_1 x_3, x_3^3\}$ pointwise converges to an Hermite projector $P$
with
\begin{align*}
\mathrm{ran} P'=\{&\delta_{(0,0,0)}\circ \mathrm{D}^{(0,0,0)}, \delta_{(0,0,0)}\circ \mathrm{D}^{(1,0,0)}, \delta_{(0,0,0)}\circ \mathrm{D}^{(0,1,0)}, \delta_{(0,0,0)} \circ \mathrm{D}^{(0,0,1)},\\
&\delta_{(1,1,1)}\circ \mathrm{D}^{(0,0,0)}, \delta_{(1,1,1)}\circ \mathrm{D}^{(1,0,0)},
\delta_{(1,1,1)}\circ \mathrm{D}^{(0,1,0)}, \delta_{(1,1,1)} \circ \mathrm{D}^{(0,0,1)}\},
\end{align*}
as $h$ tends to zero.
\end{exmp}



Finally, we  select test functions
\begin{align*}
f_1(x_1,x_2)&=1+(1-x_1)^4+(1-x_2)^4,\\
f_2(x_1,x_2,x_3)&=1+(1-x_1)^2+(1-x_2)^2+(1-x_3)^2,\\
\end{align*}
to illustrate the  pointwise convergence of ideal projectors of type
partial derivative in the above examples.

For Example 2, when
$h=1/10, 1/100, 1/1000, \ldots$, we have
\begin{align*}
P_{\frac{1}{10}}f_1=&3-\frac{385039}{99000}x_2-\frac{3439}{1000}x_1+\frac{719}{150}x_2^2+ \frac{86}{25} x_1 x_2-\frac{1438}{495}x_2^3,\\
P_{\frac{1}{100}}f_1=&3-\frac{39984109399}{9999000000}x_2-\frac{3940399}{1000000}x_1+\frac{970199}{165000}x_2^2+\frac{9851}{2500} x_1 x_2\\
&-\frac{970199}{249975}x_2^3,\\
P_{\frac{1}{1000}}f_1=&3-\frac{571426287284857}{142857000000000}x_2-\frac{3994003999}{1000000000}x_1+\frac{997001999}{166500000}
x_2^2\\
&+\frac{998501}{250000} x_1 x_2-\frac{142428857}{35714250}x_2^3,\\
\cdots&\\
P f_1=&3-4 x_2-4 x_1+6 x_2^2+ 4 x_1 x_2-4 x_2^3.
\end{align*}

For Example \ref{ex2}, 
\begin{align*}
P_{\frac{1}{10}}f_2=&4-\frac{829}{495}x_3-\frac{19}{10}x_2-\frac{19}{10}x_1-\frac{7}{3}x_3^2+2 x_2 x_3+ 2 x_1 x_3+\frac{80}{99}x_3^3,\\
P_{\frac{1}{100}}f_3=&4-\frac{989299}{499950}x_3-\frac{199}{100}x_2-\frac{199}{100}x_1-\frac{37}{33} x_3^2+2 x_2 x_3+ 2x_1 x_3\\
&+\frac{800}{9999}x_3^3,\\
P_{\frac{1}{1000}}f_2=&4-\frac{998992999}{499999500}x_3-\frac{1999}{1000}x_2-\frac{1999}{1000}x_1-\frac{337}{333}x_3^2+2
x_2 x_3+ 2 x_1 x_3\\
&+\frac{8000}{999999}x_3^3,\\
\cdots&\\
P f_2=&4-2 x_3-2 x_2-2 x_1-x_3^2+2 x_2 x_3+ 2 x_1 x_3.
\end{align*}


\section{Proof of Theorem \ref{mianthm}}\label{proof}

First of all, we need to relate forward differences of multivariate polynomials to their partial derivatives. The following formula is quite useful for this purpose.

\begin{lem}\label{L1}

Let $i, m\in \mathbb{N}_0$ satisfying $i \geq m >0$. Then
\begin{equation}\label{L1g}
\sum\limits_{j=0}^{i-1}(-1)^j {i \choose j}(i-j)^m=\left\{
                                                \begin{array}{ll}
                                                  i!, &m=i; \\
                                                  0, & m<i.
                                                \end{array}
                                              \right.
\end{equation}
\end{lem}

\begin{pf}
The proof can be completed by induction on $m$.\qed

\end{pf}


\begin{lem}\label{uni}
Let $\xi, h\in \mathbb{R}, h\neq 0, $ and $i, \alpha\in \mathbb{N}_0$. Then for every monomial $x^\alpha$ in $\mathbb{R}[x]$,
\begin{equation}\label{unishi}
\sum\limits_{j=0}^i (-1)^j {i \choose j} \delta_{\xi+h(i-j)}x^\alpha=\left\{
                                                            \begin{array}{ll}
                                                           h^{i}\delta_{\xi}\circ\mathrm{D}^ix^\alpha, & \alpha\leq i; \\
                                             h^{i}\delta_{\xi}\circ\mathrm{D}^ix^\alpha+O(h^{i+1}), &\alpha
                                             > i,
                                                            \end{array}
                                                          \right.
\end{equation}
where the remainder $O(h^{i+1})$  is a polynomial in $h$.
\end{lem}

\begin{pf}
From the theory of finite difference(see for example \cite{Ame1977}) we know that
$$\Delta^i \delta_\xi f(x)=\sum\limits_{j=0}^i (-1)^j {i \choose j} \delta_{\xi+h(i-j)}f(x)=
h^{i}\delta_{\xi}\circ \mathrm{D}^if(x)+O(h^{i+1}),$$
where $\Delta$ is the forward difference operator and $f(x)\in C^i(\mathbb{R})$. When $f(x)$ is substituted by $x^\alpha$ in this equation, (\ref{unishi}) follows immediately. Moreover, by Lemma \ref{L1}, we can
easily check that the remainder $O(h^{i+1})$ in (\ref{unishi}) is a
polynomial in $h$.  This completes the proof. \qed

\end{pf}

The conclusion of Lemma \ref{uni} will be carried over to multivariate cases as follows.

\begin{lem}\label{polyn}
Suppose that $h\in\mathbb{R}\backslash\{0\}$,
$\bm{\xi}=(\xi_1,\ldots,\xi_d) \in \mathbb{R}^d$, and
$\bm{i}=(i_1,\ldots,i_d) \in \mathbb{N}_0^d$. Then for artitrary
monomial $\bm{x}^{\bm{\alpha}}$ in $\mathbb{R}[\bm{x}]$, we have
\begin{eqnarray}\label{polyngongshi}
\sum\limits_{\bm{0}\leq\bm{j}\leq \bm{i}} (-1)^{\bm{j}}
{\bm{i}\choose \bm{j}} \delta_{\bm{\xi}+h
(\bm{i}-\bm{j})}{\bm{x}}^{\bm{\alpha}} =\left\{\begin{array}{ll}
h^{\|\bm{i}\|_1}
\delta_{\bm{\xi}}\circ\mathrm{D}^{\bm{i}}\bm{x}^{\bm{\alpha}}+O(h^{\|\bm{i}\|_1+1}),
&
\bm{i}< \bm{\alpha};\\
h^{\|\bm{i}\|_1}\delta_{\bm{\xi}}\circ\mathrm{D}^{\bm{i}}\bm{x}^{\bm{\alpha}}, &
\mbox{otherwise},
\end{array}
                                                          \right.
\end{eqnarray}
where $(-1)^{\bm{j}}=(-1)^{j_1}\cdots(-1)^{j_d}$ and ${\bm{i}\choose \bm{j}}={i_1\choose j_1}\cdots {i_d\choose j_d}$ provided that $\bm{j}=(j_1, \ldots, j_d)$.
\end{lem}

\begin{pf}
First, it follows from Lemma \ref{uni} that for every $1\leq
k\leq d$
\begin{equation}\label{danbianyuan}
\sum\limits_{j_k=0}^{i_k} (-1)^{j_k} {i_k \choose j_k}\delta_{\xi_k+h(i_k-j_k)}x_k^{\alpha_k}
=\left\{
                                                            \begin{array}{ll}
                                                           h^{i_k}\delta_{\xi_k}\circ \mathrm{D}^{i_k}x_k^{\alpha_k}, & \alpha_k\leq i_k; \\
           h^{i_k}\delta_{\xi_k}\circ \mathrm{D}^{i_k}x_k^{\alpha_k}+O(h^{i_k+1}),
&\alpha_k>i_k.
                                                            \end{array}
                                                         \right.
\end{equation}
Further, we observe that
\begin{equation}\label{hebin2}
\prod\limits_{k=1}^d\delta_{\xi_k}\circ \mathrm{D}^{i_k}x_k^{\alpha_k}=\delta_{\bm{\xi}}\circ\mathrm{D}^{\bm{i}}\bm{x}^{\bm{\alpha}}
\end{equation}
and
\begin{equation}\label{hebin}
\sum\limits_{\bm{0}\leq\bm{j}\leq \bm{i}} (-1)^{\bm{j}}
{\bm{i}\choose \bm{j}} \delta_{\bm{\xi}+h
(\bm{i}-\bm{j})}{\bm{x}}^{\bm{\alpha}}
=\prod_{k=1}^d\left(\sum\limits_{j_k=0}^{i_k}(-1)^{j_k} {i_k\choose j_k}
\delta_{\xi_k+h (i_k-j_k)}{x_k^{\alpha_k}}\right).
\end{equation}
Finally, we distinguish three cases to prove that the right-hand sides of (\ref{hebin}) and (\ref{polyngongshi}) are equal to each other, which will complete the proof.

Case 1: $\bm{\alpha}\leq \bm{i}$.

Using (\ref{danbianyuan}) and (\ref{hebin2}), it is straightforward to
verify that
$$\prod_{k=1}^d\left(\sum\limits_{j_k=0}^{i_k}(-1)^{j_k}
{i_k\choose j_k} \delta_{\xi_k+h (i_k-j_k)}{x_k^{\alpha_k}}\right)=\prod_{k=1}^d
h^{i_k}\delta_{\xi_k}\circ \mathrm{D}^{i_k}x_k^{\alpha_k}=h^{\|\bm{i}\|_1} \delta_{\bm{\xi}}\circ\mathrm{D}^{\bm{i}}\bm{x}^{\bm{\alpha}}.$$

Case 2:  $\bm{i}\not<\bm{\alpha}$ and
$\bm{\alpha}\not\leq \bm{i}$.

In this case, there must exist some
$1\leq k, l\leq d$ such that $\alpha_k< i_k$ and
$i_l<\alpha_l$. Thus, it is easily checked that
$$\prod_{k=1}^d\left(\sum\limits_{j_k=0}^{i_k}(-1)^{j_k}
{i_k\choose j_k} \delta_{\xi_k+h (i_k-j_k)}{x_k^{\alpha_k}}\right)=h^{\|\bm{i}\|_1} \delta_{\bm{\xi}}\circ \mathrm{D}^{\bm{i}}\bm{x}^{\bm{\alpha}}=0.$$

Case 3: $\bm{i}< \bm{\alpha}$.

Let $l=\max\{k: i_k< \alpha_k, 1\leq k\leq d\}$. Then, applying (\ref{danbianyuan}) and
(\ref{hebin2}), we deduce that
\begin{align*}
&\prod_{k=1}^d\left(\sum\limits_{j_k=0}^{i_k}(-1)^{j_k} {i_k\choose j_k}
\delta_{\xi_k+h (i_k-j_k)}{x_k^{\alpha_k}}\right)\\
=&\prod_{k=1}^l\left(\sum\limits_{j_k=0}^{i_k}(-1)^{j_k}
{i_k\choose j_k} \delta_{\xi_k+h (i_k-j_k)}{x_k^{\alpha_k}}\right)
\prod_{k=l+1}^d\left(\sum\limits_{j_k=0}^{i_k}(-1)^{j_k}
{i_k\choose j_k} \delta_{\xi_k+h (i_k-j_k)}{x_k^{\alpha_k}}\right)\\
=&\prod_{k=1}^l \left( h^{i_k}\delta_{\xi_k}\circ
\mathrm{D}^{i_k}x_k^{\alpha_k}+O(h^{i_k+1})\right)\prod_{k=l+1}^d
h^{i_k}\delta_{\xi_k}\circ \mathrm{D}^{i_k}x_k^{\alpha_k}\\
=&h^{\|\bm{i}\|_1}
\delta_{\bm{\xi}}\circ\mathrm{D}^{\bm{i}}\bm{x}^{\bm{\alpha}}+O(h^{\|\bm{i}\|_1+1}),
\end{align*}
where the empty product is understood to be 1. \qed
\end{pf}

Equation (\ref{polyngongshi}) makes a connection between
the forward difference calculus and the differential calculus for multivariate monomials. From Lemma \ref{uni}, it follows that the remainder $O(h^{\|\bm{i}\|_1+1})$ in (\ref{polyngongshi}) is a polynomial in $h$. Equipped with these facts, we can establish the relationship between forward differences and partial derivatives of
multivariate polynomials, which plays an important
role in the proof of Theorem \ref{mianthm}.

\begin{cor}\label{polynomialcor}
Let $\bm{i}, h, \bm{\xi}$ be as in \emph{Lemma \ref{polyn}} and $p\in \mathbb{R}[\bm{x}]\backslash \{0\}$. Then
\begin{eqnarray}\label{henxinc}
\frac{1}{ h^{\|\bm{i}\|_1}}\sum\limits_{\bm{0}\leq
\bm{j}\leq\bm{i}} (-1)^{\bm{j}}
{\bm{i}\choose \bm{j}}\delta_{\bm{\xi}+h (\bm{i}-\bm{j})} p
=\left\{
   \begin{array}{ll}
    \delta_{\bm{\xi}}\circ\mathrm{D}^{\bm{i}}p +O(h), &  p<_m\bm{i};\\
      \delta_{\bm{\xi}}\circ\mathrm{D}^{\bm{i}}p , & \hbox{otherwise}.
   \end{array}
 \right.
\end{eqnarray}
\end{cor}
\begin{pf}
Assume that nonzero polynomial $p$ has form (\ref{pform}). Since
$$\sum\limits_{\bm{0}\leq \bm{j}\leq\bm{i}} (-1)^{\bm{j}}
{\bm{i}\choose \bm{j}}\delta_{\bm{\xi}+h (\bm{i}-\bm{j})} p
=\sum\limits_{\bm{\alpha}}\widehat{p}({{\bm{\alpha}}})\sum\limits_{\bm{0}\leq
\bm{j}\leq\bm{i}} (-1)^{\bm{j}}
{\bm{i}\choose \bm{j}}\delta_{\bm{\xi}+h (\bm{i}-\bm{j})}
{\bm{x}}^{{\bm{\alpha}}}$$ and
$$\delta_{\bm{\xi}}\circ\mathrm{D}^{\bm{i}}p=\sum\limits_{\bm{\alpha}}\widehat{p}({{\bm{\alpha}}})
\delta_{\bm{\xi}}\circ\mathrm{D}^{\bm{i}} {\bm{x}}^{{\bm{\alpha}}},$$
we get
\begin{align*}
\sum\limits_{\bm{0}\leq \bm{j}\leq\bm{i}} (-1)^{\bm{j}}
{\bm{i}\choose \bm{j}}\delta_{\bm{\xi}+h (\bm{i}-\bm{j})} p
=\left\{
   \begin{array}{ll}
   { h^{\|\bm{i}\|_1}} \delta_{\bm{\xi}}\circ\mathrm{D}^{\bm{i}}p +O( h^{\|\bm{i}\|_1+1}), &  p <_m \bm{i};\\
      { h^{\|\bm{i}\|_1}}\delta_{\bm{\xi}}\circ\mathrm{D}^{\bm{i}}p , & \hbox{otherwise},
   \end{array}
 \right.
\end{align*}
which leads to the corollary immediately.\qed
\end{pf}

Now, we are ready to prove Theorem \ref{mianthm}.

\vskip 8pt

\noindent\textsc{Proof of Theorem \ref{mianthm}.} We adopt the notation of Definition \ref{HermiteProjectorde} and Remark \ref{dualbasis}. Let $\bm{q}=(q_1,q_2,\ldots,q_n)$ be an $\mathbb{R}$-basis for $\mathrm{ran}P$. Without loss of generality, we assume that the entries of $\bm{\lambda}$ and
$\bm{\lambda}_h$ are ordered ascendingly w.r.t. $\prec_\lambda$ and then denoted as $o_1, \ldots, o_n$ and $o'_1, \ldots, o'_n$ respectively. For convenience, we set $n \times n$ matrices
$$
\bm{\lambda}^T \bm{q}=(o_i q_j)_{1\leq i, j\leq n}, \quad \bm{\lambda}^T_h \bm{q}=(o'_i q_j)_{1\leq i, j\leq n},
$$
and, therefore, for every $q \in \mathbb{R}[\bm{x}]$, $n$ by $1$ vectors
$$
\bm{\lambda}^T q=(o_i q)_{1\leq i\leq n}, \quad \bm{\lambda}^T_h q=(o'_i q)_{1\leq i\leq n}.
$$
By Corollary \ref{polynomialcor}, equation (\ref{henxinc}) can be rewritten as
$$
\delta_{\bm{\xi}}\circ\mathrm{D}^{\bm{i}}p=\left\{
                                       \begin{array}{ll}
                                         \frac{1}{ h^{\|\bm{i}\|_1}}\sum\limits_{\bm{0}\leq \bm{j}\leq\bm{i}}
(-1)^{\bm{j}} {\bm{i}\choose \bm{j}}\delta_{\bm{\xi}+h (\bm{i}-\bm{j})}
p+O(h), &  p <_m \bm{i}; \\
                                         \frac{1}{ h^{\|\bm{i}\|_1}}\sum\limits_{\bm{0}\leq \bm{j}\leq\bm{i}}
(-1)^{\bm{j}} {\bm{i}\choose \bm{j}}\delta_{\bm{\xi}+h (\bm{i}-\bm{j})}
p, & \hbox{otherwise,}
                                       \end{array}
                                     \right.
$$
which implies that for fixed $1\leq k \leq \mu$ and $\bm{i} \in
\mathfrak{\Delta}^{{(k)}}$, $\delta_{\bm{\xi}^{(k)}}\circ\mathrm{D}^{\bm{i}}p$ can be
linearly expressed by $\{\delta_{\bm{\xi}^{(k)}+h \bm{l}} p: \bm{l}\in \mathfrak{\Delta}^{(k)}\}\cup \{O(h)\}$ since $\mathfrak{\Delta}^{{(k)}}$ is lower, and moreover, the linear combination
coefficient of each $ \delta_{\bm{\xi}^{(k)}+h \bm{l}}p$ is independent of $p\in \mathbb{R}[\bm{x}]$. Thus, it turns out that there exists a nonsingular matrix $T_p$ of order $n$ such that
\begin{equation}\label{juzhen}
\left[\widehat{\bm{\lambda}_h^T \bm{q}}\Big| \widehat{\bm{\lambda}_h^T
q}\right]:=T_p\left[\bm{\lambda}^T_h \bm{q}| \bm{\lambda}^T_h
q\right]=\left[\bm{\lambda}^T \bm{q}| \bm{\lambda}^T q\right]+\left[E_h|
\bm{\epsilon}_{h}\right],
\end{equation}
where each entry of  $[E_h|\bm{\epsilon}_{h}]$ is either $0$ or
$O(h)$. As a consequence, the linear systems
$$\left(\widehat{\bm{\lambda}_h^T \bm{q}}\right) \bm{x}=\widehat{\bm{\lambda}_h^T
q}\quad \mbox{and} \quad\left(\bm{\lambda}^T_h \bm{q}\right) \bm{x}=\bm{\lambda}^T_h q$$
are equivalent, namely they have the same set of solutions.

(i) From (\ref{juzhen}), it follows that each entry of
matrix $\widehat{\bm{\lambda}_h^T \bm{q}}$ converges to its
corresponding entry of matrix $\bm{\lambda}^T \bm{q}$  as $h$ tends
to zero, which implies that
$$\lim\limits_{h\rightarrow
0}\det\left(\widehat{\bm{\lambda}_h^T \bm{q}}\right)=\det \left(\bm{\lambda}^T
\bm{q}\right).$$
Since $\det (\bm{\lambda}^T \bm{q})\neq 0$, there exists $\eta>0$ such that
$$\det\left(\widehat{\bm{\lambda}_h^T \bm{q}}\right) \neq 0, \quad0<|h|<\eta.$$
Notice that (\ref{juzhen}) directly leads to  $\mathrm{rank}
\left(\widehat{\bm{\lambda}_h^T \bm{q}}\right)=\mathrm {rank}\left(\bm{\lambda}_h^T
\bm{q}\right)$,
$$
\mathrm{ran}P_h=\mathrm{Span}_\mathbb{R}\bm{q}, \quad 0<|h|<\eta,
$$
follows, i.e., $\bm{q}$ forms an $\mathbb{R}$-basis for $\mathrm{ran}P_h$. Since $\bm{q}$ is also a basis for $\mathrm{ran}P$, we have
$$
\mathrm{ran}P=\mathrm{ran}P_h,\quad 0<|h|<\eta.
$$

(ii) Suppose that $\widetilde{\bm{x}}_h$ and $\widetilde{\bm{x}}$ be the unique solutions of nonsingular linear systems
\begin{equation}\label{Leq}
(\bm{\lambda}_h^T \bm{q})\bm{x}=\bm{\lambda}_h^T q
\end{equation}
and
\begin{equation}\label{Heq}
(\bm{\lambda}^T \bm{q})\bm{x}=\bm{\lambda}^T q
\end{equation}
respectively, where $0<|h|<\eta$. It is easy to see that
$$
P_h q=\bm{q\widetilde{\bm{x}}}_h\quad \mbox{and}\quad Pq=\bm{q\widetilde{\bm{x}}}.
$$
Remark that, as $h\rightarrow 0$,  $P$ is the pointwise limit of $P_h$ if and only if $P q$ is the coefficientwise limit of $P_h q$ for all $q\in \mathbb{R}[\bm{x}]$. Therefore, it is sufficient to show that for every $q\in \mathbb{R}[\bm{x}]$, the solution vector of system (\ref{Leq}) converges to the one of system (\ref{Heq}) when $h$ tends to zero, namely
$$
\lim_{h\rightarrow 0} \bm{\widetilde{x}}_h= \bm{\widetilde{x}}.
$$
By (\ref{juzhen}), the linear system
\begin{equation}\label{widehat}
    \left(\widehat{\bm{\lambda}_h^T
\bm{q}}\right)\bm{x}=\widehat{\bm{\lambda}_h^T q}
\end{equation}
can be rewritten as
$$\left(\bm{\lambda}^T \bm{q}+E_h\right)\bm{x}=\left(\bm{\lambda}^T
q+\bm{\epsilon}_{h}\right).
$$
Since system (\ref{widehat}) is equivalent to system (\ref{Leq}), $\widetilde{\bm{x}}_h$ is also the unique solution of it. Consequently, applying the perturbation analysis of the sensitivity
of linear systems (see for example \cite{Matrixcompu1996}, p.80ff), we have
$$\left\|\widetilde{\bm{x}}_h-\bm{\widetilde{x}}\right\|\leq\left\|{(\bm{\lambda}^T \bm{q})}^{-1}\right\|\left\|\bm{\epsilon}_h-E_h\bm{\widetilde{x}}\right\|+O(h^2).$$
Since each entry of vector $\bm{\epsilon}_h-E_h\bm{\widetilde{x}}$  is
either $0$ or $O(h)$, it follows that $\lim\limits_{h\rightarrow
0}\|\widetilde{\bm{x}}_h-\bm{\widetilde{x}}\|= 0$, or, equivalently,
$\lim\limits_{h\rightarrow 0}\widetilde{\bm{x}}_h=\bm{\widetilde{x}}$, which completes the proof of the theorem. \qed

\bibliographystyle{elsarticle-num}
\bibliography{ref}

\begin{thebibliography}{10}
\expandafter\ifx\csname url\endcsname\relax
  \def\url#1{\texttt{#1}}\fi
\expandafter\ifx\csname urlprefix\endcsname\relax\def\urlprefix{URL }\fi
\expandafter\ifx\csname href\endcsname\relax
  \def\href#1#2{#2} \def\path#1{#1}\fi

\bibitem{dBo94}
C.~de~Boor, Polynomial interpolation in several variables, in: Studies in
  Computer Science, Plenum Press, New York, 1994, pp. 87--119.

\bibitem{GS00:2}
M.~Gasca, T.~Sauer, Polynomial interpolation in several variables, Adv. Comput.
  Math. 12 (2000) 377--410.

\bibitem{She2009}
B.~Shekhtman, Ideal interpolation: Translations to and from {Algebraic}
  {Geometry}, in: L.~Robbiano, J.~Abbott (Eds.), Approximate Commutative
  Algebra, Texts and Monographs in Symbolic Computation, Springer Vienna, New
  York, 2009, pp. 163--192.

\bibitem{Bir1979}
G.~Birkhoff, The algebra of multivariate interpolation, in: C.~V. Coffman,
  G.~J. Fix (Eds.), Constructive Approaches to Mathematical Models, Academic
  Press, New York, 1979, pp. 345--363.

\bibitem{dBo2005}
C.~de~Boor, Ideal interpolation, in: C.~K. Chui, M.~Neamtu, L.~L. Schumaker
  (Eds.), Approximation Theory XI: Gatlinburg 2004, Nashboro Press, Brentwood,
  TN, 2005, pp. 59--91.

\bibitem{LZD2011}
Z.~Li, S.~Zhang, T.~Dong, On the existence of certain error formulas for a
  special class of ideal projectors, J. Approx. Theory 169~(9) (2011)
  1080--1090.

\bibitem{BS2006}
B.~Shekhtman, On a conjecture of {Carl} de {Boor} regarding the limits of
  {Lagrange} interpolants, Constr. Approx. 24~(3) (2006) 365--370.

\bibitem{Foga1968}
J.~Fogarty, Algebraic families on an algebraic surface, Amer. J. Math. 90
  (1968) 511--521.

\bibitem{deBoorShe2008}
C.~de~Boor, B.~Shekhtman, On the pointwise limits of bivariate {Lagrange}
  projectors, Linear Algebra Appl. 429~(1) (2008) 311--325.

\bibitem{BS2008}
B.~Shekhtman, Bivariate ideal projectors and their perturbations, Adv. Comput.
  Math. 29~(3) (2008) 207--228.

\bibitem{Lor1992}
R.~Lorentz, Multivariate Birkhoff Interpolation, Vol. 1516 of Lecture Notes in
  Mathematics, Springer, Heidelberg, 1992.

\bibitem{deboor2006}
C.~de~Boor, What are the limits of {Lagrange} projectors?, in: B.~D. Bojanov
  (Ed.), Constructive Theory of Functions (Varna 2005), Marin Drinov Acad.
  Publ. House, Sofia, 2006, pp. 51--63.

\bibitem{Mor2009}
T.~Mora, Gr\"{o}bner technology, in: M.~Sala, T.~Mora, L.~Perret, S.~Sakata,
  C.~Traverso (Eds.), Gr\"{o}bner Bases, Coding, and Cryptography, Springer,
  Berlin, 2009, pp. 11--25.

\bibitem{CM1995}
L.~Cerlienco, M.~Mureddu, From algebraic sets to monomial linear bases by means
  of combinatorial algorithms, Discrete Math. 139~(1-3) (1995) 73--87.

\bibitem{FRR2006}
B.~Felszeghy, B.~R\'{a}th, L.~R\'{o}nyai, The lex game and some applications,
  J. Symbolic Comput. 41~(6) (2006) 663--681.

\bibitem{Lor2000}
R.~Lorentz, Multivariate {Hermite} interpolation by algebraic polynomials: A
  survey, J. Comput. Appl. Math. 122 (2000) 166--201.

\bibitem{Ame1977}
W.~F. Ames, Numerical Methods for Partial Differential Equations, Academic
  Press, New York, 1977.

\bibitem{Matrixcompu1996}
G.~H. Golub, C.~V. Loan, Matrix Computations, third ed., The Johns Hopkins
  University Press, Baltimore, MD, 1996.

\end{thebibliography}

\end{document}